\documentclass[12pt]{amsart}

\usepackage[colorlinks=true, pdfstartview=FitV, linkcolor=blue,
citecolor=blue, urlcolor=blue]{hyperref}

\usepackage{mathtools}
\usepackage{amsfonts}
\usepackage{amsmath}
\usepackage{amssymb} 
\usepackage{amsthm}

\usepackage{times}
\usepackage[T1]{fontenc}
\usepackage{enumitem}
\usepackage{setspace}
\usepackage{microtype}
\usepackage{cite}

\allowdisplaybreaks

\usepackage[margin=1.2in]{geometry}

\newtheorem{theorem}{Theorem}[section]

\theoremstyle{definition}

\theoremstyle{remark}
\newtheorem{remark}[theorem]{Remark}

\numberwithin{equation}{section}

\newcommand{\R}{\mathbb R}
\newcommand{\N}{\mathbb N}

\newcommand{\rn}{{{\mathbb R}^n}}
\newcommand{\rd}{{{\mathbb R}^d}}

\DeclareMathOperator{\supp}{supp}

\DeclareMathOperator*{\esssup}{ess\,sup}

\newcommand{\D}{\mathcal D}
\newcommand{\Ss}{\mathcal S}

\newcommand{\Rdf}{\mathcal{R}}

\DeclareMathOperator{\op}{op}

\newcommand{\A}{\mathcal A}
\newcommand{\K}{\mathcal K}
\newcommand{\W}{\mathcal W}

\DeclareMathOperator{\clconv}{\overline{conv}}

\newcommand{\vecf}{\vec{f}}

\newcommand{\Sd}{\mathcal{S}_d}

\def\Xint#1{\mathchoice
   {\XXint\displaystyle\textstyle{#1}}%
   {\XXint\textstyle\scriptstyle{#1}}%
   {\XXint\scriptstyle\scriptscriptstyle{#1}}%
   {\XXint\scriptscriptstyle\scriptscriptstyle{#1}}%
   \!\int}
\def\XXint#1#2#3{{\setbox0=\hbox{$#1{#2#3}{\int}$}
     \vcenter{\hbox{$#2#3$}}\kern-.5\wd0}}

\def\avgint{\Xint-}

\begin{document}

\title[Factorization and Extrapolation]
{Matrix weights, singular integrals, \\Jones factorization and \\
    Rubio de Francia extrapolation}

\author{David Cruz-Uribe, OFS}
\address{Dept. of Mathematics \\
University of Alabama \\
 Tuscaloosa, AL 35487, USA}
\email{dcruzuribe@ua.edu}


\date{March 20, 2024}

\thanks{The author is partially supported by a Simons Foundation
  Travel Support for Mathematicians  Grant.  He would like to thank
  Yoshihiro Sawano for his invitation and generous support to speak at
  the Conference on Research on Real, Complex and Functional Analysis
  in Kyoto, Japan, October 26--27, 2023.  The original article was
  based on his talks at the conference; this revised version is
  updated to include very recent work on the $\A_2$ conjecture.}

\begin{abstract}
  In this article we give an overview of the problem of finding sharp
  constants in matrix weighted norm inequalities for singular
  integrals, the so-called matrix $A_2$ conjecture.  We begin by
  reviewing the history of the problem in the scalar case, including a
  sketch of the proof of the scalar $A_2$ conjecture.  We then discuss
  the original, qualitative results for singular integrals with matrix
  weights and the best known quantitative estimates.  We give an
  overview of new results by the author and Bownik, who developed a
  theory of harmonic analysis on convex set-valued functions.  This
  led to the proof the Jones factorization theorem and the Rubio de
  Francia extrapolation theorem for matrix weights, two longstanding
  problems.  Rubio de Francia extrapolation was expected to be a major
  tool in the proof of the matrix $\A_2$ conjecture; however, this conjecture was very
  recently proved false.  We discuss this problem.
\end{abstract}

\maketitle

\section{Introduction}
One weight norm inequalities have been extensively studied since the
work of Muckenhoupt and others in the
1970s.  (See~\cite{duoandikoetxea01,garcia-cuerva-rubiodefrancia85,grafakos08b}
for details and references.)  Central to their study is the
Muckenhoupt $A_p$ condition.  It gives a necessary and sufficient
condition for many of the operators of classical harmonic analysis
(maximal operators, singular integral operators, square functions) to
satisfy weighted norm inequalities.   There are two very important results in
the study of of $A_p$ weights.  The first is  the Jones factorization
theorem~\cite{Jones}, which gives a complete characterization of the
class $A_p$ in terms of the simpler class $A_1$.  The second is the Rubio de
Francia extrapolation theorem~\cite{Rubio}, which, in its
simplest form, shows that if an operator is bounded on weighted $L^2$, then
it is bounded on weighted $L^p$ for all $p$, $1<p<\infty$.   Rubio de
Francia's colleague Antonio Cordoba~\cite{garcia-cuerva87} summarized this result by
remarking that it showed that  {\em there are no $L^p$ spaces, only weighted $L^2$.}

Since the 1990s, there has been a great deal of interest in extending
the theory of scalar weights to the setting of matrix weights:  that
is, $d\times d$ measurable matrix functions $W$ that are self-adjoint,
positive semi-definite, and act on vector-valued functions $\vecf$.
These problems were first studied by Nazarov, Treil and
Volberg~\cite{MR1428988,MR1428818,MR1423034,MR1478786}, who asked if the
weighted norm inequalities for singular integrals could be extended to
the matrix setting.  They defined a class of matrix weights, $\A_p$, and
showed that the Hilbert transform is bounded with respect to this
class.  This result was later extended to general singular integral
operators by Christ and Goldberg~\cite{MR2015733,CG}.

More recently, attention has been focused on determining the sharp
constant in matrix norm inequalities.  In the scalar case,
Hyt\"onen~\cite{hytonenP2010} proved that the sharp constant in the
weighted $L^p$ norm inequality is proportional to
$[w]_{A_p}^{\max\{1,p'-1\}}$.  In the scalar case, most proofs of this
result proceed in two steps: first, it is proved in the case $p=2$,
and then it is extended to all $p$ by using a sharp version of the
Rubio de Francia extrapolation theorem.
(See~\cite{MR3204859,Lern2012,MR2140200}.) For this reason, the
problem is referred to as the ``$A_2$ conjecture.''

For more than a
decade it has been an open
question as to whether the same bound is true for matrix weights.  The
best known result is that it is bounded above by
$[W]_{\A_p}^{1+\frac{1}{p-1}-\frac{1}{p}}$, which, when $p=2$ gives an
exponent of $\frac{3}{2}$ rather than the conjectured value of $1$.
Very recently, Domelevo, Petermichl, Treil and
Volberg~\cite{DPTV-2024} show that this conjecture is false and that
  when $p=2$, $\frac{3}{2}$ is the best possible exponent.  The
  question of sharp exponents when $p\neq 2$ is still open.

In the 1990s
Nazarov and Treil~\cite{MR1428988} posed two fundamental problems:
extend the Jones factorization theorem and the Rubio de Francia
extrapolation theorem to matrix weights.  These problems were only
solved in 2022 by the author and Bownik~\cite{mb-2022}.    While no
longer directly applicable to the matrix $\A_2$ conjecture, this
result is a powerful tool for proving matrix weighted inequalities
when $p\neq 2$.  

\medskip

The goal of this article is to provide an overview of the theory of
matrix weights and of the proofs of
factorization and extrapolation for matrix weights.   It is organized
as follows.  In Section~\ref{section:scalar} we provide some
background about scalar weights.  We focus on those results which are
essential for understanding the subsequent development of the matrix
theory.  In Section~\ref{section:matrix} we review some of the history
of matrix weights.  We will concentrate on the results of Christ and
Goldberg~\cite{MR2015733,CG} and the more recent work of Nazarov,
Petermichl, Treil, and Volberg~\cite{MR3689742}.   In
Section~\ref{section:convex} we will discuss the work of the author
and Bownik.  Our approach, which considerably extends the work of
Nazarov, {\em et al.}, builds a theory of harmonic analysis on
convex set-valued functions. This in turn provides the necessary tools to
prove factorization and extrapolation for matrix weights.  Finally, in
Section~\ref{section:guess} we make some final remarks about
exrapolation and the matrix $A_2$ conjecture.

Throughout this paper we will use the following notation.  In
Euclidean space the constant $n$ will denote the dimension of $\R^n$,
which will be the domain of our functions.  The value $d$ will denote
the dimension of vector and set-valued functions.  For
$1\leq p \leq \infty$, $L^p(\R^n)$ will denote the Lebesgue space of
scalar functions, and $L^p(\R^n,\rd )$ will denote the Lebesgue space
of vector-valued functions.

By a cube we will always mean a cube in $\rn$ with sides parallel to
the coordinate axes.  The Lebesgue measure of a cube, or of any
arbitrary set $E$, will be denoted
by $|E|$.  By a weight we mean a non-negative, locally
integrable function that is positive except on a set of measure $0$.
We define $w(E) = \int_E w(x)\,dx$, and we let $\avgint_Q w(x)\,dx =
|Q|^{-1}w(Q)$.  By $L^p(w)$ we mean the scalar weighted space with
measure $w\,dx$, $L^p(\rn,w\,dx)$. 

Given $v=(v_1,\ldots,v_d)^t \in \rd $, the Euclidean norm of $v$ will
be denoted by $|v|$; from context there should be no confusion with
the notation for Lebesgue measure of a set.  The closed unit ball in
$\{ v \in \rd : |v| \leq 1\}$ will be denoted by ${\mathbf B}$.
Matrices will be $d\times d$ matrices with real-valued entries.  The
set of all $d\times d$, symmetric, positive semidefinite matrices will
be denoted by $\Sd$.

Given two quantities $A$ and $B$, we will write
$A \lesssim B$, or $B\gtrsim A$ if there is a constant $c>0$ such that
$A\leq cB$.  If $A\lesssim B$ and $B\lesssim A$, we will write
$A\approx B$.

 \section{Background: the scalar theory of weights}
 \label{section:scalar}

 In harmonic analysis, a fundamental class of weights are those that
 satisfy the Muckenhoupt $A_p$ condition, introduced
 in~\cite{muckenhoupt72}.  Given $1<p<\infty$, $w\in A_p$ if
\[ [w]_{A_p} = \sup_Q \left(\avgint_Q w(x)\,dx\right) \left(\avgint_Q
  w(x)^{1-p'}\,dx\right)^{p-1} < \infty, \]
where the supremum is taken over all cubes.  A weight $w$ is in $A_1$ if
\[ [w]_{A_1} = \sup_{Q} \esssup_{Q\ni x}  w(x)^{-1} \avgint_Q
w(y)\,dy < \infty. \]

The quantity $[w]_{A_p}$ is referred to as the $A_p$ characteristic of
a weight.  These weights arise naturally in the study of Hardy-Littlewood maximal
operator,
\[ Mf(x) = \sup_Q \avgint_Q |f(y)|\,dy \cdot \chi_Q(x), \]
where the supremum is taken over all cubes.  For $1\leq p<\infty$, the
$A_p$ condition is necessary and sufficient for the maximal operator
to satisfy the weak $(p,p)$ inequality
\[ w(\{ x\in \rn : Mf(x) > \lambda \})
  \leq
  \frac{C}{t^p} \int_\rn |f(x)|^p w(x)\,dx, \]
and for $1<p<\infty$ it is necessary and sufficient for the strong
$(p,p)$ inequality
\[ \int_\rn Mf(x)^p w(x)\,dx
  \leq
  C \int_\rn |f(x)|^p w(x)\,dx. \]

  Recall that a Calder\'on-Zygmund singular integral $T$
is a bounded
operator on $L^2(\rn)$ for which there exists a kernel $K(x,y)$, defined on
$\rn\times \rn\setminus \Delta$, where $\Delta =\{ (x,x) : x \in \rn\}$, such that
if $f \in L_c^\infty(\rn)$ and $x\not\in \supp(f)$, then
\[ Tf(x) = \int_\rn K(x,y)f(y)\,dy.  \]
The kernel satisfies the size and regularity conditions
\[ |K(x,y)| \leq \frac{C}{|x-y|^n}, \]
\[ |K(x+h,y)-K(x,y)|+|K(x,y+h)-K(x,y)|
  \leq
  C\frac{|h|^\delta}{|x-y|^{n+\delta}}, \]
where $|x-y|>2|h|$.  The $A_p$ condition is also sufficient for a
singular integral operator to satisfy the weak and strong $(p,p)$
inequalities; moreover, it is necessary for non-degenerate singular
integrals such as the Riesz transforms.
See~\cite{duoandikoetxea01,garcia-cuerva-rubiodefrancia85,
  grafakos08b,stein93}.

The weak $(p,p)$ inequality for singular integrals can be proved using
kernel estimates and the good/bad decomposition of Calder\'on and
Zygmund (see, for instance,~\cite{garcia-cuerva-rubiodefrancia85}).
The strong-type inequality was originally proved by comparing the norm
of the singular integral operator to that of the maximal operator.
Coifman and Fefferman~\cite{coifman-fefferman74} proved that given
$p$, $0<p<\infty$, and $w\in A_q$ for any $q$, $1\leq q<\infty$, there
exists a constant depending on $[w]_{A_q}$ such that
\begin{equation} \label{eqn:TM}
\int_\rn T^*f(x)^p w(x)\,dx
  \leq
  C\int_\rn Mf(x)^p w(x)\,dx.  
\end{equation}
Here, $T^*$ is the maximal singular integral, defined by
\[ T^* f(x) = \sup_{\epsilon>0}|T_\epsilon f(x)|
  = \sup_{\epsilon>0}
  \bigg|\int_{|x-y|>\epsilon} K(x,y)f(y)\,dy\bigg|, \]
which dominates the singular integral pointwise.  They proved this by
proving a so-called good-$\lambda$ inequality:  there exists
$\delta>0$ such that for every $\gamma,\,\lambda>0$ and for every cube $Q$,
\[ w(\{ x\in Q : T^*f(x) >2\lambda, Mf(x)\leq \gamma\lambda\})
  \leq
  C\gamma^\delta w(Q). \]
An alternative proof of \eqref{eqn:TM} using the sharp maximal
operator was given by
Journ\'e~\cite{MR706075}; see also Alvarez and
P\'erez~\cite{alvarez-perez94}.

For the past three decades there has been a great deal of interest in
determining the best constant in the strong $(p,p)$ inequality for
singular integrals in terms of the $A_p$
characteristic of $w$.  This question was first considered by Buckley in the 1990s~\cite{buckley93}; it became the subject of concerted effort when Astala,
Iwaniec and Saksman~\cite{MR1815249} proved that sharp regularity
results for solutions of the Beltrami equation hold provided that the
Beurling-Ahlfors operator satisfies $\|Tf\|_{L^p(w)} \leq C
[w]_{A_p}\|f\|_{L^p(w)}$ for $p\geq 2$.  This problem was extended to
all Calder\'on-Zygmund operators and all $p>1$:  it was conjectured
that 
\begin{equation} \label{eqn:sharp} 
\|Tf\|_{L^p(w)} \leq C(n,T,p) [w]_{A_p}^{\max(1, p'-1)}\|f\|_{L^p(w)}.
\end{equation}
By the Rubio de Francia extrapolation theorem, which we will discuss
below, it suffices to prove this when $p=2$; for this reason this
problem was referred to as the $A_2$ conjecture.
 It was studied by a number of authors, including Lacey, 
Petermichl, and Volberg~\cite{lacey-petermichl-reguera2010,
 MR2354322,petermichl08, petermichl-volberg02}, and was finally solved
in 2010 by Hyt\"onen~\cite{hytonenP2010,hytonen-perez-treil-volbergP}.
His proof was quite difficult and we will not consider it. 

Lerner and Nazarov~\cite{lerner-IMRN2012,MR3127380,MR4007575} and
Conde-Alonso and Rey~\cite{MR3521084} gave 
new and simpler proofs of the $A_2$ conjecture.  As part of their proofs
they introduced the technique of sparse domination.  Let
$\D$ be any translation of the standard dyadic grid.  A collection $\Ss\subset
\D$ is said to be sparse if there exists a collection of pairwise
disjoint sets $\{ E(Q) : Q \in \Ss \}$ such that for each $Q$, $E(Q)
\subset Q$ and $|Q|\leq 2|E(Q)|$.   A sparse operator is an averaging
operator of the form
\[ T_{\Ss}f(x) = \sum_{Q\in \Ss} \left( \avgint_Q
  f(y)\,dy\right)\cdot \chi_Q(x). \]
They showed that given a bounded function of compact support, there
exist a finite collection of dyadic grids $\{\D_n\}_{n=1}^N$ and
sparse families $\{\Ss_n\}_{n=1}^N$ such that
\[  |Tf(x)| \leq C \sum_{n=1}^N T_{\Ss_n}(|f|)(x).  \]
Given this, the $A_2$ conjecture reduces to proving the corresponding
estimates for sparse operators.  This is done in two steps.  First,
weighted $L^2$ bounds are proved using an argument
from~\cite{dcu-martell-perez}.  Let $w\in A_2$ and let $\sigma=w^{-1}$.
By duality, there exists $h\in L^2(w)$, $\|h\|_{L^2(w)}=1$, such that
    \begin{align*}
       \|T_\Ss f\|_{L^2(w)}
       & = \int_{\rn} T_\Ss f(x) h(x) w(x)\,dx \\
   & \leq 2 \sum_{Q\in \Ss} \avgint_Q f(x)\,dx
       \avgint_Q h(x)w(x)\,dx |E(Q)| \\
  & = 2 \sum_{Q\in \Ss} \frac{w(Q)}{|Q|}\frac{\sigma(Q)}{|Q|}
     \frac{1}{\sigma(Q)}\int_Q f(x)w(x)\sigma(x)\,dx \frac{1}{w(Q)} \int_Q h(x)w(x)\,dx |E(Q)| \\
& \leq 2[w]_{A_2} \sum_{Q\in \Ss} \int_{E(Q)}
 \tilde{M}_\sigma^d (fw)(x) \sigma(x) \tilde{M}_w^dh(x) w(x) \,dx \\
& \leq 2[w]_{A_2} \int_\rn \tilde{M}_\sigma^d (fw)(x) \sigma(x) \tilde{M}_w^dh(x) w(x) \,dx \\
& \leq 2[w]_{A_2} \|\tilde{M}_\sigma^d (fw)\|_{L^2(\sigma)}
                 \|\tilde{M}_w^d h\|_{L^2(w)} \\ 
   & \leq 8[w]_{A_2} \|fw\|_{L^2(\sigma)}
        \| h\|_{L^2(w)} \\
 & = 8[w]_{A_2} \|f\|_{L^2(w)}. 
     \end{align*}
Here, $\tilde{M}_w^d$ is the ``universal'' dyadic maximal operator defined with respect to
the grid $\D$ and the measure $w\,dx$:
\[ \tilde{M}_w^d f(x) =\sup_{Q\in \D} \frac{1}{w(Q)}\int_Q
  |f(y)|w(y)\,dy \cdot \chi_Q(x),  \]
where the supremum is taken over all dyadic cubes in the grid $\D$.  This
operator is bounded on $L^2(w)$ with a constant independent of $w$
(see~\cite{garcia-cuerva-rubiodefrancia85}).  Similarly, $\tilde{M}_\sigma^d$
is bounded on $L^2(\sigma)$.

The proof that sparse operators satisfy weighted $L^p$ bounds with the
desired constant follows at once from the Rubio de Francia
extrapolation theorem.  Here we state it in a
``sharp constant'' version first proved by~ Dragi\v{c}evi\'{c}, {\em et
  al.}~\cite{MR2140200}.  (See also~\cite{cruz-martell-perezBook}.)

\begin{theorem} \label{thm:rubio}
Given $p_0$, $1\leq p_0<\infty$, suppose that for some operator $T$
and 
for all $w_0\in A_{p_0}$, the inequality
\[ \int_\rn |Tf(x)|^{p_0}w_0(x)\,dx
  \leq
  N_{p_0}([w_0]_{A_{p_0}})\int_\rn |f(x)|^{p_0}
w_0(x)\,dx \]
holds.  Then for all $p$, $1<p<\infty$,  and all $w\in A_p$, 
\[ \int_\rn |Tf(x)|^{p}w(x)\,dx
  \leq
   N_{p}([w]_{A_{p}})\int_\rn |f(x)|^{p}
w(x)\,dx. \]
\end{theorem}

The proof of extrapolation ultimately depends on three things:
\begin{enumerate}

  \item The duality of $A_p$ weights:  for $1<p<\infty$, $w\in A_p$ if
    and only if $\sigma=w^{1-p'} \in A_{p'}$, and
    $[\sigma]_{A_{p'}}=[w]_{A_p}^{p'-1}$.  This is an immediate
    consequence of the definition of $A_p$ weights.

  \item A quantitative bound for the maximal operator: for
    $1<p<\infty$,
    $ \|Mf\|_{L^p(w)} \leq C[w]_{A_p}^{p'-1}\|f\|_{L^p(w)}$.  This
    was first explicitly proved by Buckley~\cite{buckley93} but was
    implicit in Christ and Fefferman~\cite{MR684636}.

 \item The Jones factorization theorem:  for $1<p<\infty$, $w\in A_p$
   if and only if there exist $w_0,\,w_1 \in A_1$ such that
   $w=w_0w_1^{1-p}$.   This was first proved by Jones~\cite{Jones},
   but a much more elementary proof was given by Coifman, Jones, and
   Rubio de Francia~\cite{MR687639}.  (See also~\cite{dcu-paseky}.)

 \end{enumerate}
 
With the first two facts we can define the Rubio de Francia iteration
operator, which is also fundamental to the proof of the Jones
factorization theorem.  
Let $M$ be  the Hardy-Littlewood maximal operator and $w\in A_p$,
$1<p<\infty$.  Given a
non-negative function $h$,  define
\[ \Rdf h(x) = \sum_{k=0}^\infty \frac{M^kh(x)}{2^k
    \|M\|_{L^p(w)}^k}. \]
It then follows from the definition and the properties of $A_p$ weights that
\begin{enumerate}

\item $h(x) \leq \Rdf h(x)$;

\item $\|\Rdf h\|_{L^p(w)} \leq 2 \| h\|_{L^p(w)}$;

  \item $\Rdf h \in A_1$ and $[\Rdf h]_{A_1} \leq 2\|M\|_{L^p(w)} \leq
    C[w]_{A_p}^{p'-1}$.

  \end{enumerate}

  We briefly sketch the proof of extrapolation.  For simplicity,
  we will only consider the case when $p_0=2$, and we will not give
  the proof that yields the best possible constant.  (For complete
  details, plus references to other proofs,
  see~\cite{dcu-paseky,cruz-martell-perezBook}.)  Fix $p$, $p\neq 2$,
  and $w\in A_p$.   Let $\Rdf_1$ be the Rubio de Francia iteration
  operator defined above, and let $\Rdf_2$ be the operator
  corresponding to $\sigma=w^{1-p'} \in A_{p'}$.  Then by duality
  there exists $h\in L^{p'}(w)$, $\|h\|_{L^{p'}(w)}=1$, such that 
\begin{align*}
  \|Tf\|_{L^p(w)}
  & = \int_\rn |Tf(x)| h(x)w(x)\,dx \\
& \leq \int_\rn |Tf(x)| 
\Rdf_1f(x)^{-\frac{1}{2}} \Rdf_1f(x)^{\frac{1}{2}}   \Rdf_2 (hw)(x)\,dx \\
& \leq \left(\int_\rn |Tf(x)|^2 \Rdf_1f(x)^{-1} \Rdf_2
  (hw)(x)\,dx\right)^{\frac{1}{2}} 
\left(\int_\rn \Rdf_1f(x) \Rdf_2 (hw)(x)\,dx \right) ^{\frac{1}{2}} \\
& = I_1 ^{\frac{1}{2}} I_2 ^{\frac{1}{2}}. 
\end{align*}
We estimate each term separately.  The estimate of $I_2$ uses the properties of the
iteration operators:
\begin{align*}
  I_2
  & = \int_\rn \Rdf_1f(x) w(x)^{\frac{1}{p}}\Rdf_2 (hw)(x)
      w(x)^{-\frac{1}{p}}\,dx \\ 
& \leq \|\Rdf_1f\|_{L^p(w)} 
\|\Rdf_2(hw)\|_{L^{p'}(\sigma)} \\
& \leq 4\|f\|_{L^p(w)} \|hw\|_{L^{p'}(\sigma)}\\
& \leq 4\|f\|_{L^p(w)}.
\end{align*}
To estimate $I_1$, we use our hypothesis, the fact that by the
Jones factorization theorem, $\Rdf f_1(f)^{-1} \Rdf f_2\in A_2$,  and the
estimate for $I_2$:
\begin{align*}
  I_1
  & = \int_\rn |Tf(x)|^2 \Rdf_1f(x)^{-1} \Rdf_2
  (hw)(x)\,dx \\
& \leq C \int_\rn |f|^2 \Rdf_1f(x)^{-1} \Rdf_2
  (hw)(x)\,dx \\
& \leq C \int_\rn |f(x)|\Rdf_2
  (hw)(x)\,dx \\
& \leq CI_2.
\end{align*}
If we combine these inequalities, we get the desired result.

\begin{remark}
  Note that in this proof of extrapolation, we only use the easier
  half of the Jones factorization theorem:  that if $w_0,\,w_1\in
  A_1$, then $w_0w_1^{1-p} \in A_p$.  This property, referred to as
  ``reverse factorization'' (see~\cite{cruz-martell-perezBook}),
  follows at once from the definition of the $A_p$ and $A_1$
  conditions.  We will return to this fact below.
\end{remark}

\section{Matrix weights and matrix-weighted norm inequalities}
\label{section:matrix}

In the 1990s, Nazarov, Treil and Volberg in a series of
papers~\cite{MR1428988,MR1428818,MR1478786,MR1423034} considered the
question of whether the theory of Muckenhoupt weights could be
extended to matrix weights applied to vector-valued functions.  Beyond
the intrinsic interest of this problem, their original motivation came
from problems in the study of multivariate random stationary
processes, and from the study of Toeplitz operators acting on
vector-valued Hardy spaces.

To describe the problem, we define some notation.  Let
$\vecf=(f_1,\ldots,f_d)$.  Given a singular integral $T$, define  it
acting on a vector-valued function by
\[ T\vecf = (Tf_1,\ldots,Tf_d).  \]
It is immediate that if $1<p<\infty$ and $\vecf \in L^p(\rn,\rd )$,
then $\|T\vecf\|_{ L^p(\rn,\rd )} \leq C \|\vecf\|_{ L^p(\rn,\rd )}$.

To define weights, recall that $\Sd$ is the set of $d\times d$,
self-adjoint, positive semi-definite matrices.  A matrix weight is a
measurable function  $W : \rn \rightarrow \Sd$.  We define an
associated scalar weight using the operator norm on $W$:
\[ |W(x)|_{\op}
  =
  \sup_{\substack{\xi \in \rd \\|\xi|=1}} |W(x)\xi|. \]
We define the matrix weighted space $L^p(W)=L^p(W,\rn,\rd )$ with the norm
\[ \|\vecf\|_{L^p(W)} = \bigg( \int_\rn
  |W(x)^{\frac1p}\vecf(x)|^p\,dx\bigg)^{\frac1p}. \]
Note that when $d=1$, this reduces to the scalar space $L^p(w)$.  Also
note that $W^{\frac1p}$ is well-defined since $W$ is positive
semi-definite. (We note in passing that in~\cite{mb-2022} we defined
this space with $W^{\frac1p}$ replaced by $W$.   We believe that this is the
correct way to define matrix weighted spaces, but here, for consistency with
the earlier literature, we use the historical definition.)

With this notation, the problem posed by Nazarov, Treil and Volberg is
the following:  given $1<p<\infty$, prove there is a Muckenhoupt-type
condition on matrix weights so that the inequality
\[ \int_{\R^n} |W^{\frac{1}{p}}(x)T{f}(x)|^p \,dx
  \leq C \int_{\R^n} |W^{\frac{1}{p}}(x){f} (x)|^p \,dx \]
holds for singular integrals.    Treil and Volberg~\cite{MR1428818}  first
solved this problem on the real line for the Hilbert transform  when
$p=2$.  They showed that this inequality holds if $W$ satisfies an analog of the $A_2$
condition:
\begin{equation} \label{eqn:matrix-A2}
[W]_{\A_2} =
  \sup_Q \bigg| \left(\avgint_Q W(x)\,dx\right)^{\frac{1}{2}}
  \left(\avgint_Q W^{-1}(x)\,dx\right)^{\frac{1}{2}}\bigg|_{\op}
  <\infty. 
\end{equation}
This condition, however, does not extend to the case $p\neq 2$.  An
equivalent, but more technical definition of matrix $\A_p$ in terms of
norm functions was
conjectured by Treil and used by Nazarov and Treil~\cite{MR1428988}
and by Volberg~\cite{MR1423034} to prove matrix-weighted norm
inequalities for the Hilbert transform.  Their idea was to replace
matrices with norms on $\rd $.  Here we sketch their definition; for
complete details see the above references or~\cite{mb-2022}.

Let  $\rho : \rn \times \rd 
\rightarrow [0,\infty)$ be a measurable function such that for
a.e. $x\in \rn$, $\rho(x,\cdot)$ is a norm on $\rd $:  that is, given
$v,\,w\in \rd$ and $\alpha \in \R$,
\begin{enumerate}
\item $\rho(x,v)=0$ if and only if $v=0$;
\item $\rho(x,v+w) \leq \rho(x,v)+\rho(x,w)$;
\item $\rho(x,\alpha v)=|\alpha|\rho(x,v)$. 
\end{enumerate}
For instance, given a matrix weight $W$, we can define a norm by
$\rho_W(x,v)=|W(x)v|$. 
Define the dual of a norm function $\rho^*$ by
\[ \rho^*(x,v) = \sup_{w\in \rd , \rho(x,w)\leq 1} |\langle v,
  w\rangle|. \]
Finally, given a cube $Q$, and $1\leq p<\infty$, define the average of a
norm function on $Q$ by
\[ \langle \rho \rangle_{p,Q}(v)
  =
  \bigg(\avgint_Q \rho(x,v)^p\,dx\bigg)^{\frac1p}. \]
With this notation, they defined a norm $\rho$ to be in $\A_p$ if for
every $v\in \rd $,
\[  \langle \rho^* \rangle_{p',Q}(v)
  \leq
  C  \langle \rho \rangle_{p,Q}(v)^* \]
When $d=1$ this immediately reduces to the Muckenhoupt $A_p$
condition.  

Later,
Roudenko~\cite{MR1928089} gave an equivalent definition of 
$\A_p$ that more closely resembled the scalar definition:  for
$1<p<\infty$, $W\in \A_p$ if 
\begin{equation} \label{eqn:roudenko1}
 [W]_{\A_p} = \sup_Q \avgint_Q \bigg( \avgint_Q
    |W^{\frac{1}{p}}(x)W^{-\frac{1}{p}}(y)|_{\op}^{p'}\,dy\bigg)^{\frac{p}{p'}}\,dx <
    \infty. 
  \end{equation}
  Frazier and Roudenko~\cite{MR2104276} also introduced the concept of matrix $\A_1$
  weights:  $W\in \A_1$ if
  \begin{equation} \label{eqn:roudenko2}
  [W]_{\A_1} = \sup_Q \esssup_{x\in Q}
    \avgint_Q |W^{-1}(x)W(y)|_{\op} \,dy < \infty. 
  \end{equation}
  
Norm inequalities for   Calder\'on-Zygmund singular
integrals in $\R^n$ were proved by Christ and
Goldberg~\cite{CG,MR2015733}.  Treil and Volberg had earlier noted that a key
obstacle to proving matrix-weighted inequalities was the lack of a
maximal operator that did not lose the geometric
information embedded in a vector-valued function.  A key
component of the proofs in~\cite{CG,MR2015733} is a scalar-valued, matrix
weighted maximal operator (now referred to as the Christ-Goldberg
maximal operator):
\[ M_W{\vecf}(x) =
  \sup_Q \avgint_Q |W(x)W^{-1}(y){\vecf}(y)|\,dy\cdot \chi_Q(x).  \]
The motivation for this definition comes from the following
observation:  if $T$ is a singular integral operator (or indeed any
linear operator), then $T$ satisfies $\|T\vecf\|_{L^p(W)} \leq
C\|\vecf\|_{L^p(W)}$ if and only if the operator $T_W$, defined by
\[ T_W\vecf(x) = W(x)^{\frac1p}T(W^{-\frac1p}\vecf)(x) \]
satisfies $\|T\vecf\|_{L^p(\rn,\rd )} \leq
C\|\vecf\|_{L^p(\rn,\rd )}$.

The first step was to prove that if $W\in \A_p$, $1<p<\infty$, then $M_W$ is bounded
from  $L^p(\rn,\rd)$ into $L^p(\rn)$.  This required the introduction
of an auxiliary maximal operator $M_W'$, defined by
\[ M_W{\vecf}(x) =
  \sup_Q \avgint_Q |\W_Q^p W^{-1}(y){\vecf}(y)|\,dy\cdot \chi_Q(x).  \]
Here, $\W_Q^p$ is a constant matrix defined as follows:  given the
norm $\langle \rho_W \rangle_{p,Q}$, let $K_W^p$ be the closed unit ball
in $\rd$ with respect to this norm:
\[ K_W^p = \{ v \in \rd : \langle \rho_W \rangle_{p,Q}(v) \leq 1 \}. \]
Then $K_W^p$ is a convex set in $\rd$, so there exists a unique ellipsoid
of maximal volume, called the John ellipsoid and that we denote by
$E_Q^p$, such that $E_Q^p \subset K_W^p \subset \sqrt{d} E_Q^p$.
Finally, there exists a matrix $\W_Q^p$ such that
$E_Q^p=\W_Q^p{\mathbf B}$, and the norm induced by this matrix is
equivalent to $\langle \rho_W \rangle_{p,Q}$. 
 The matrix $\W_Q^p$ is
referred to as the reducing matrix associated to $W$ on $Q$.
Note that when $d=1$, $\W_Q^p$  is just the $p$-average of the
weight.

They showed that the auxiliary maximal operator maps $L^p(\rn,\rd)$ into
$L^p(\rn)$; this can be done using Calder\'on-Zygmund cubes and an argument analogous to
that for the maximal operator.  Then, via a  stopping time
argument, they showed that
$\|M_W\vecf\|_{L^p(\rn,\rd)} \leq C \|M_W\vecf\|_{L^p(\rn,\rd)}$.

Given this maximal operator, to prove weighted norm inequalities for
singular integrals they adapted the ideas of Coifman and Fefferman to
prove a good-$\lambda$ inequality.  Define
\[ T_W^*\vecf(x) = \sup_{\epsilon>0}
  |W^{\frac1p}(x)T_\epsilon(W^{-\frac1p}\vecf)(x)|. \]
Then for every smooth function $\vecf$ with compact support, they
proved that there exist constants $0<b<1$ and $c>0$ such that for all
$\lambda>0$,
\begin{multline*}
  |\{x\in \rn : T_W^*\vecf(x)>\lambda,
  \max\{ M_W'\vecf(x), M_W\vecf(x)\}<c\lambda\}| \\
  \leq
  \frac{1}{2}b^p|\{x\in \rn : T_W^*\vecf(x) > b\lambda\}|. 
\end{multline*}

\medskip

As in the scalar case, this approach to proving matrix weighted norm
inequalities does not yield quantitative estimates on the constant in
terms of the $[W]_{\A_p}$ characteristic.  After the proof of the
  scalar $A_2$ conjecture by Hyt\"onen,  it was conjectured that the
  corresponding result holds for matrix weights:  for $1<p<\infty$, 
  \[ \|T\vecf\|_{L^p(W)}
    \leq
    C[W]_{\A_p}^{\max\{1,p'-1\}}\|\vecf\|_{L^p(W)}.  \]
  This problem was first considered by Bickel, Petermichl and
Wick~\cite{MR3452715} and by Pott and Stoica~\cite{MR3698161} when
$p=2$.  More recently,  Nazarov, Petermichl, Treil and
Volberg~\cite{MR3689742} showed that when $p=2$, they can get a
constant of the form $C(n,d,T)[W]_{A_2}^{\frac32}$.  (Also
see~\cite{MR3818613}.)  Their proof is based on a deep
generalization of the sparse domination estimates described above in
the scalar case; they then estimate their sparse operator using square
function estimates.    

The sparse operator introduced in~\cite{MR3689742} replaces
vector-valued functions by averages that are convex sets.  They show
that given $\vecf \in L_c^\infty(\rn,\rd)$, there exists a finite
collection of sparse sets  $\{\Ss_n\}_{n=1}^N$ such that
\begin{equation} \label{eqn:nptv1}
 T{\vecf}(x) \in C\sum_{n=1}^N\sum_{Q\in \Ss_n}
  \langle\langle{\vecf}\rangle\rangle_Q \chi_Q(x), 
\end{equation}
where $\langle\langle{f}\rangle\rangle_Q$ is the convex set
\[  \langle\langle{\vecf}\rangle\rangle_Q
  = \bigg\{ \avgint_Q k(y) {\vecf}(y)\,dy : k \in L^\infty(Q),
  \|k\|_\infty \leq 1 \bigg\},  \]
and the sum is the (infinite) Minkowski
sum of convex sets (see below for a definition).   They referred to this estimate as a convex body
sparse domination.  To complete their proof, however, they did not
work directly with these convex set-valued functions, but rather
replaced them by vector-valued sparse operators of the form
\[ T^S {\vecf}(x) = \sum_{Q\in \Ss} \avgint_Q
  \varphi_Q(x,y){\vecf}(y)\,dy, \]
where for each $Q$, $\varphi_Q$ is a real-valued function supported on
$Q\times Q$  such that, for each $x$,
$\|\varphi_Q(x,\cdot)\|_\infty \leq 1$.
These they estimated using square function estimates.

In the case $p\neq 2$,  quantitative
results were proved by  the author,
Isralowitz and Moen~\cite{MR3803292}, who got a constant of the form
\[ C(n,d,p,T)[W]_{A_p}^{1+\frac{1}{p-1}-\frac{1}{p}}.  \]
They used the sparse domination result of Nazarov, {\em et al.} and
reduced to the vector-valued sparse operators, but
instead of using square function estimates, they adapted techniques from the
theory of $A_p$ bump conditions in the study of scalar, two-weight
norm inequalities.  In both proofs there appeared to be a loss of information:
even in the scalar case these techniques do not yield the $A_2$
conjecture.  So the problem becomes to work directly with the convex
set-valued sparse operator.

\section{Convex set-valued functions and matrix weights}
\label{section:convex}

The convex body sparse operator of
Nazarov {\em et al.} introduced a powerful new tool into the study of
matrix weights.  Building upon these ideas, we further developed the
theory of convex set-valued functions to prove Jones factorization
and Rubo de Francia extrapolation for matrix weights.   This problem
was first raised by Nazarov and Treil~\cite{MR1428988}:
\begin{quote}
    {\em Actually, the whole theory of scalar $(A_p)$-weights  can be
      transferred to the matrix case except two results...the Peter
      Jones factorization and the Rubio-de-Francia extrapolation
      theory.  But today (May 1, 1996) we do not know what the
      analogues of these two things are in high dimensions.}
    \end{quote}
    Their assessment at that time was somewhat optimistic: for
    instance, there is still not a complete theory of reverse H\"older
    inequalities for matrix weights (though
    see~\cite{davey2022exponential}).  Nevertheless, they did identify
    two fundamental problems in the study of matrix weights.  While
    factorization and extrapolation are interesting in and of
    themselves, they gained greater importance with the study of the
    matrix $\A_2$ conjecture.  As we described above, extrapolation
    played an important role in proving the scalar $A_2$ conjecture.
    Moreover, for a number of technical reasons matrix $\A_2$ weights
    are easier to work with than matrix $\A_p$ weights, $p\neq
    2$. (E.g., the simpler definition~\eqref{eqn:matrix-A2} can be
    used.)  Therefore, it seems natural to try to prove a sharp
    version of extrapolation in the matrix case.

    Both results have been fully proved for matrix
    weights.  Here we state the two main theorems from~\cite{mb-2022}.

    \begin{theorem} \label{thm:matrix-jones}
Given $1<p<\infty$, then $W\in \A_p$ if and only if
          there exist commuting matrices $W_0,\,W_1\in \A_1$ such that
         \[ W = W_0W_1^{1-p}. \]
     \end{theorem}

\begin{remark}
  For simplicity and ease of comparison to the scalar case, we state
  Theorem~\ref{thm:matrix-jones} assuming that the
  matrices $W_0$ and $W_1$ 
  commute.  We can remove this hypothesis, but to do so we must
  replace the product $W_0W_1^{1-p}$ with a more complicated
  expression, the geometric  mean of the two matrices:  see
  \cite{bhatia}. 
\end{remark}

\medskip

\begin{theorem} \label{thm:matrix-rubio}
Given $p_0$, $1\leq p_0<\infty$, suppose that for some operator $T$ and every $W_0\in
      \A_{p_0}$,
      \[ \bigg(\int_{\mathbb R^n} |W_0^{\frac{1}{p_0}}
        T\vecf|^{p_0} \,dx)\bigg)^{\frac{1}{p_0}}
        \leq N_{p_0}([W_0]_{p_0})\bigg( \int_{\mathbb R^n}
        |W_0^{\frac{1}{p_0}} \vecf|^{p_0}
        \,dx \bigg)^{\frac{1}{p_0}}.\]
      Then for all $p$, $1<p<\infty$, and $W\in \A_p$,
      \[ \bigg(\int_{\mathbb R^n}
        |W^{\frac{1}{p}} T\vecf|^{p}\,dx)\bigg)^{\frac{1}{p}}
        \leq N_{p}([W]_{p})\bigg( \int_{\mathbb R^n}
        |W^{\frac{1}{p}} \vecf|^{p}
        \,dx \bigg)^{\frac{1}{p}}.\]
\end{theorem}

\smallskip

\begin{remark}
  We actually prove a more general version of
  Theorem~\ref{thm:matrix-rubio}, replacing the operator
  $T$ by a family of pairs of functions $(f,g)$.  This more abstract
  approach to extrapolation was first suggested
  in~\cite{cruz-uribe-perez00} and systematically developed
  in~\cite{cruz-martell-perezBook}.
\end{remark}

\begin{remark}
  In Theorem~\ref{thm:matrix-rubio} the function $N_p$
  depending on $N_{p_0}$ has exactly the same form as the functions
  gotten in the sharp constant extrapolation theorem of
  Dragi{\v{c}}evi{\'c}, {\em et al.}~\cite{MR2140200}.  
\end{remark}

\medskip

The proofs of Theorems~\ref{thm:matrix-rubio}
and~\ref{thm:matrix-jones} are both long and extremely technical, and
it is beyond the scope of this article to give many details.  Instead
we will provide a conceptual overview of the proofs and of the tools
we developed for them.  As part of the proof we have begun to lay the
groundwork for studying harmonic analysis on convex set-valued
functions.  These results are of interest in their own right.
Moreover, beyond their application to proving these two theorems, we
believe that they will be useful for exploring more deeply the convex
body sparse bounds of Nazarov, {\em et al.}.

Underlying our proofs of Theorems~\ref{thm:matrix-jones}
and~\ref{thm:matrix-rubio} was the systematic philosophy of trying to
replicate the proofs of factorization and extrapolation in the scalar
case, particularly the elementary proofs of these results
in~\cite{cruz-martell-perezBook,dcu-paseky}.  As we noted above, the
proofs of extrapolation and factorization depended on
the duality of scalar  $A_p$ and sharp bounds for weighted norm
inequalities for the maximal operator; these in turn were used to
build the Rubio de Francia iteration operator.   Extrapolation further
required the Jones factorization theorem, though as we saw in the
proof sketched above, we only used the easier direction of this
result.

The fundamental technical obstacle to the proof was the lack of an
appropriate definition of the maximal operator.  While the
Christ-Goldberg maximal operator $M_W$ was sufficient to prove strong
$(p,p)$ bounds for singular integrals, it has the drawback that it
maps a vector-valued function $\vecf$ to scalar-valued function
$M_W\vecf$.  Therefore, it cannot be iterated, and so cannot be used
to construct a Rubio de Francia iteration operator.  To overcome this
problem, we passed from vector-valued functions to the larger category
of convex set-valued functions, and defined a convex set-valued
maximal operator.

We begin with some definitions related to convex sets.  For complete
details, see~\cite{mb-2022} and the references it contains.  Let $\K$
denote the family of all convex sets $K\subset \rd$ that are closed,
bounded, and symmetric: i.e., if $x\in K$, then $-x\in K$.  Sometimes
it is necessary also to assume $K$ is absorbing: that is, that $0\in
\mathrm{int}(K)$.  However, we will not worry about this technical
hypothesis.  Given a set $K\subset \rd$, let $|K|=\{ |v| : v \in K
\}$; given a matrix $W$, define $WK= \{ Wv : v \in K\}$; note that
$WK$ is also a convex set.   Given two convex sets $K$ and $L$, their
Minkowski sum is the set $K+L = \{ u+ v : u \in K, v\in L \}$. 

As we noted above, every norm has associated to it the convex set
which is its unit ball.  The converse is also true:  to  every convex
set $K$ these is associated to it a unique
norm, $\rho_K$ on $\rd$.  Moreover, arguing as we did before, there
exists a matrix $W$ such that $\rho_K \approx
\rho_W$.  While we generally want to work with matrices $W\in \A_p$,
at key points in the argument it is necessary to pass to working with
more general norms and the underlying convex sets.

We now define convex set-valued functions $F$; there is actually a
well-developed theory of such objects, but it does not appear to be
well-known among harmonic analysts:  see~\cite{MR2458436,CV}.  Let $F: \rn
\rightarrow \K$ be such a map.  There are several equivalent ways to define
measurability of such functions $F$; for our purposes a useful and
intuitive definition is that there exists a countable family of
functions $\vecf_k : \rn \rightarrow\rd$ such that for almost every $x$,
\[ F(x) = \overline{ \{ \vecf_k(x) : k \in \N \} }. \]
Such functions are referred to as selection functions.

Given such a convex set-valued function $F$, we can define for each
$x$ the associated John Ellipsoid $E(x)$.  This ellipsoid-valued
function is measurable. This fact has been used in the literature
(see, for instance, Goldberg~\cite{MR2015733}) but we could not find a
proof in the literature; a proof is given in~\cite{mb-2022}.

The integral of a convex set-valued function can be defined using
selection functions--this object is referred to as the Aumann
integral~\cite{MR0185073} (see also~\cite{MR2458436}).  Given $\Omega
\subset \rn$ and a function $F : \Omega \rightarrow \K$, define
\[ S^1(\Omega, F) = \{ \vecf \in L^1(\Omega,\rd) : \vecf(x) \in
  F(x)\}.  \]
Then the Aumann integral of $F$ is defined to be the set
\[ \int_\Omega F(x)\,dx
  =
  \bigg\{ \int_\Omega \vecf(x)\,dx : \vecf \in S^1(\Omega,F)
  \bigg\}. \]
It can be shown that since $F(x)$ is closed, bounded and convex, the
Aumann integral is also a closed, convex set in $\rd$. 

There is a close connection between the Aumann integral and the convex
averages $\langle\langle \vecf \rangle\rangle_Q$ used by Nazarov, {\em
  et al.} to define their convex body sparse operator.  Given a
vector-valued function $\vecf$, define $F_{\vecf}(x)$ to be the closed
convex hull of the set $\{\vecf(x), -\vecf(x)\}$.  Then $F_{\vecf}$ is a
measurable, convex set-valued function, and
\[ \langle\langle \vecf \rangle\rangle_Q
  = \avgint_Q F_{\vecf}(x)\,dx. \]
Thus, the convex body sparse operator $T_\Ss$ is a convex set-valued
function as defined above.   

We use the Aumann integral to define a convex set-valued maximal
operator.  Given $F : \rn \rightarrow \K$, let
\[ MF(x) = \clconv\bigg(
  \bigcup_Q \avgint_Q F(y)\,dy \cdot \chi_Q (x) \bigg);  \]
that is, $MF(x)$ is the closed, convex hull of the union of the Aumann integral averages of
$F$ over all cubes containing $x$.   Then $MF$ is a measurable,
convex set-valued function.   The intuition
behind this definition is that the Hardy-Littlewood maximal operator
finds the largest average in magnitude, and so uses the supremum.  The
convex set-valued maximal operator finds the largest average in
magnitude in each direction; to preserve the information about
direction we take the union of all averages.

It should be noted that this maximal
operator does not preserve some natural subsets of $\K$.  For
instance, given a vector-valued function $\vecf$, $MF_{\vecf}$
can be an absorbing convex set that properly contains the convex set
that is the closure of $\{ M\vecf,-M\vecf\}$.  If $F$ is
ellipsoid-valued--that is, $F=W{\mathbf B}$, where $W$ is a matrix
valued function--then $MF$ need not be ellipsoid-valued.

The convex set-valued maximal operator has a number of properties that
correspond to those of the Hardy-Littlewood maximal operator: for
convex set-valued functions $F$ and $G$, almost every $x\in \rn$, and $\alpha \in
[0,\infty)$, 
\begin{enumerate}

\item $F(x) \subset MF(x)$;

  \item $M(F+G)(x) \subset MF(x) + MG(x)$, where the sum is the
    Minkowski sum;

  \item $M(\alpha F)(x) = \alpha MF(x)$.
\end{enumerate}
The maximal operator $M$ also satisfies $L^p$ norm inequalities.  For
$1\leq p<\infty$, define $L^p_\K(\rn,|\cdot|)$ to be the collection of
all convex set-valued functions $F$ such that
\[ \|F\|_{L^p_\K(\rn,|\cdot|)}
  =
  \bigg(\int_\rn |F(x)|^p\,dx\bigg)^{\frac1p} < \infty.  \]
Then for $1<p<\infty$ we have that $\|MF\|_{L^p_\K(\rn,|\cdot|)} \leq
C \|F\|_{L^p_\K(\rn,|\cdot|)}$.

Weighted norm inequalities for the convex maximal operator are
governed by the matrix $\A_p$ weights.  Define
$L^p_\K(\rn,W)=L^p_\K(W)$ to be all convex set-valued functions $F$ such
that
\[ \|F\|_{L^p_\K(\rn,W)}
  =
  \bigg(\int_\rn |W^{\frac1p}(x)F(x)|^p\,dx\bigg)^{\frac1p} < \infty.  \]

\begin{theorem} \label{thm:convex-max-wts}
 Given $p$, $1<p<\infty$, and $W\in \A_p$, for every $F \in
 L^p_\K(W)$,
 \[ \|MF\|_{L^p_\K(W)}
   \leq
   C(n,d,p)[W]_{\A_p}^{p'-1}\|F\|_{L^p_\K(W)}. \]
\end{theorem}

The proof of Theorem~\ref{thm:convex-max-wts} uses the $L^p$ norm
inequalities for the Christ-Goldberg maximal operator $M_W$; the sharp
constant in terms of $[W]_{\A_p}$ was proved by Isralowitz and
Moen~\cite{MR4030471}.  It would be of interest to have a direct proof that
did not require using the Christ-Goldberg maximal operator.

We now define the analog of the $A_1$ condition for convex
set-valued functions.  Given $F : \rn \rightarrow \K$, we say that
$F\in \A_1^\K$ if $MF(x) \subset CF(x)$.  Denote the infimum of all
such constants $C$  by $[F]_{\A_1^\K}$.  There is a very close
relationship between convex set-valued $\A_1^\K$ and matrix $\A_1$.

\begin{theorem}
  Given a matrix weight $W$, $W\in \A_1$ if and only if $F=W{\mathbf
    B} \in \A_1^\K$, and $[W]_{\A_1} \approx [F]_{\A_1^\K}$. 
\end{theorem}

The proof of this result required a careful development of the
properties of norm functions and their duals on the one hand, and the
properties of convex sets and their polar bodies on the other.

With the machinery we have developed, we define the convex set-valued
analog of the Rubio de Francia iteration operator.   Given $F : \rn
\rightarrow \K$  and $W\in \A_p$, $1<p<\infty$, define
\[ \Rdf F(x) = \sum_{k=0}^\infty \frac{M^kF(x)}{2^k
    \|M\|^k_{L^p_\K(W)}}. \]
Then $\Rdf F$ has the following properties which are the exact analogs
of the properties of the scalar iteration operators:
\begin{enumerate}
\item $F(x) \subset \Rdf F(x)$;

\item $\|\Rdf F\|_{L^p_\K(W)} \leq 2 \|F\|_{L^p_\K(W)}$;

  \item $\Rdf F \in \A_1^\K$ and $[\Rdf F]_{\A_1^\K} \leq 2
    \|M\|_{L^p_\K(W)}$.
  \end{enumerate}

  \medskip

We now give an overview of the proofs of
Theorems~\ref{thm:matrix-jones} and~\ref{thm:matrix-rubio}.  The proof
of factorization in the Jones factorization theorem follows the scalar
proof in~\cite{dcu-paseky} closely.  There is one major technical
obstacle:  the scalar proof uses a variant of the Hardy-Littlewood
maximal operator, $M_sf(x) = M(|f|^s)^{\frac1s}$, $s>1$, which is a
sublinear operator.   To define $M_sF(x)$, requires replacing $F$ by an
appropriate ellipsoid (so that the power $F^s$ is defined) and then
proving that the resulting operator is sublinear.    The converse,
proving reverse factorization, is very delicate and much more
difficult than in the scalar case.  We were not able to
prove it using the definitions of matrix $\A_p$
and $\A_1$ given by Roudenko and Frazier.  For
this proof we were required to work with the definition of matrix
$\A_p$ in terms of norms as originally given by Treil and Volberg.
Intuitively, the proof can be thought of as an interpolation argument
between finite dimensional Banach spaces.

The proof of extrapolation follows the proof of sharp constant
extrapolation in the scalar case given in~\cite{cruz-martell-perezBook}.  Because of the
way in which we define the matrix $\A_p$ condition, we are also able
to include naturally in the proof  a result for extrapolation from the endpoint
$p=\infty$; this gives a quantitative version of an extrapolation result
first proved by Harboure, Mac\'\i as and Segovia~\cite{MR944321}.
Recently we learned that this quantitative version was proved earlier
by Nieraeth~\cite{nieraeth-thesis,MR4000248}.  The proof uses reverse
factorization and also requires passing between convex set-valued
functions and closely associated ellipsoid valued ones.  

\section{Extrapolation and the failure of the matrix
  $A_2$ conjecture}
\label{section:guess}

In this final section we return to the matrix $\A_2$
conjecture.  As noted above, given Theorem~\ref{thm:matrix-rubio}, it
would have been
enough to prove the $\A_2$ conjecture in the case $p=2$, as our version of matrix
extrapolation yields the sharp constants for the other values of $p$.
However, with the proof in~\cite{DPTV-2024} that this conjecture is
false when $p=2$, it is no longer clear if extrapolation can be used
to prove the sharp constant estimates for singular integrals.  If we
apply extrapolation starting with the estimate $[W]_{\A_2}^{\frac{3}{2}}$ for $p=2$. we get a
bound proportional to $[W]_{\A_p}^{\frac{3}{2}\{1,\frac{1}{p-1}\}}$,
which is worse than the estimate
$[W]_{A_p}^{1+\frac{1}{p-1}-\frac{1}{p}}$
proved in~\cite{MR3803292}.   Note that this latter exponent 
asymptotically approaches $1$ as $p\rightarrow \infty$ and
$\frac{1}{p-1}$ as $p\rightarrow 1$, the values we would have gotten
if the matrix $\A_2$ conjecture is true.  Similarly, if we begin extrapolation
by taking $p_0\neq 2$ and using the bound
$[W]_{A_{p_0}}^{1+\frac{1}{p_0-1}-\frac{1}{p_0}}$ we do not get
anything better.   It therefore remains an open question as to what
the sharp bounds are when $p\neq 2$.   It is tempting to conjecture
that the result in~\cite{MR3803292} is the best possible, but there is
no evidence to support this.  

That the exponent for singular integrals is $\frac{3}{2}$ is
suprising,  but  it is worth noting that this exponent has
  recently appeared as a lower bound in the scalar case  for the sharp
  constant for rough singular integrals when $p=2$:  see~\cite{MR4584703}.
  It is unclear if there is any connection.   In the matrix case, the
  best known exponent when $p=2$, is $\frac{5}{2}$:  see~\cite{MR4245601}.
  Starting with this, we can use our matrix extrapolation result  to get the quantitative bound
  $[W]_{\A_p}^{\frac{5}{2}\{1,\frac{1}{p-1}\}}$ for all $p$. 

  \medskip
  
Finally, we want to note that the machinery of harmonic analysis on convex set-valued
functions provides a way to work with the convex body sparse operator
directly, avoiding the loss of information in the proofs in~\cite{MR3689742,MR3803292} that
yielded the best known estimates. We believe that this approach will
be useful in other problems.  We illustrate this approach by sketching a proof  of weighted norm
inequalities for the convex
set-valued sparse operator that is quite elementary, but yields a
suboptimal constant.  By inequality~\eqref{eqn:nptv1} this yields an
elementary proof of a matrix weighted norm inequality for singular
integrals.    Given $F\in L^2_\K(W)$, by a duality argument we can
show that there exists $G \in L^2_\K(W^{-1})$, $\|G\|_{L^2_\K(W^{-1})} \leq
\sqrt{d}$, such that 

 \begin{align*}
    \|T_\Ss F\|_{L^2_\K(W)}
   & \leq \int_{\rn} \langle T_\Ss F(x), G(x)\rangle \,dx.   \\
 \intertext{If we apply the definition of $T_\Ss$, linearity, and
   Cauchy's inequality, we get}
  & \leq 2 \sum_{Q\in \Ss}
      \bigg\langle \avgint_Q F(y)\,dy, \avgint_Q G(y)\,dy
      \bigg\rangle |E(Q)| \\
   & = \sum_{Q\in \Ss}
      \int_{E(Q)} \bigg\langle  \avgint_Q
      W^{\frac12}(x)W^{-\frac12}(y) W^{\frac12}(y)F(y)\,dy, \\
& \quad       \avgint_Q
      W^{-\frac12}(x)W^{\frac12}(y) W^{-\frac12}(y)G(y)\,dy
                                                                 \bigg\rangle
                                                                 \,dx \\
    & \leq \sum_{Q\in \Ss} \int_{E_Q}  M_{W}(W^{\frac12}F)
      M_{W^{-1}}(W^{-\frac12}G)\,dx \\
   & \leq \|M_{W}(W^{\frac12}F)\|_{L^2}
     \|M_{W^{-1}}(W^{-\frac12}G)\|_{L^2},\\
   \intertext{where $M_W$ is the Christ-Goldberg maximal operator.
   By the sharp constant estimate for $M_W$,}
   & \leq C[W]_{A_2}^2 \|F\|_{L^2_\K(W)}
     \|G\|_{L^2_\K(W^{-1})} \\
   & \leq C[W]_{A_2}^2 \|F\|_{L^2_\K(W)}.
  \end{align*}

  Originally, we thought that this approach could be modified to prove
  the matrix $A_2$ conjecture.  Our goal was to modify the scalar proof sketched above, which used the
  universal dyadic maximal operator $\tilde{M}^d_w$.  We were able to  show
  that the above argument can be modified to replace the
  Christ-Goldberg maximal operators with an operator $\tilde{M}_W^d$
  that is a matrix weighted version of the universal dyadic maximal
  operator.  If this operator were bounded on $L^2(W)$ with a constant
  independent of the $\A_2$ characteristic of $W$, then we could prove
  the matrix $\A_2$ conjecture.  However, Nazarov, Petermichl,
  \v{S}kreb, and Treil~\cite{NPST-2022} showed that there
  exist matrices $W$ such that $\tilde{M}_W^d$ is not bounded on
  $L^2(W)$.   We initially beleived that this was just a technical
  obstruction, but it is clear now that this fact is fundamental.

\bibliographystyle{plain}
\bibliography{Convex-Proc ver 2}

\begin{thebibliography}{10}

\bibitem{alvarez-perez94}
J.~Alvarez and C.~P{\'e}rez.
\newblock Estimates with {$A\sb \infty$} weights for various singular integral
  operators.
\newblock {\em Boll. Un. Mat. Ital. A (7)}, 8(1):123--133, 1994.

\bibitem{MR1815249}
K.~Astala, T.~Iwaniec, and E.~Saksman.
\newblock Beltrami operators in the plane.
\newblock {\em Duke Math. J.}, 107(1):27--56, 2001.

\bibitem{MR2458436}
J.-P. Aubin and H.~Frankowska.
\newblock {\em Set-valued analysis}.
\newblock Modern Birkh\"auser Classics. Birkh\"auser Boston, Inc., Boston, MA,
  2009.
\newblock Reprint of the 1990 edition [MR1048347].

\bibitem{MR0185073}
R.~J. Aumann.
\newblock Integrals of set-valued functions.
\newblock {\em J. Math. Anal. Appl.}, 12:1--12, 1965.

\bibitem{bhatia}
R.~Bhatia.
\newblock {\em Positive definite matrices}.
\newblock Princeton Series in Applied Mathematics. Princeton University Press,
  Princeton, NJ, 2007.

\bibitem{MR3452715}
K.~Bickel, S.~Petermichl, and B.~Wick.
\newblock Bounds for the {H}ilbert transform with matrix {$A_2$} weights.
\newblock {\em J. Funct. Anal.}, 270(5):1719--1743, 2016.

\bibitem{mb-2022}
M.~Bownik and D.~Cruz-Uribe.
\newblock Extrapolation and factorization of matrix weights.
\newblock {\em preprint}, 2022.
\newblock arXiv:2210.09443.

\bibitem{buckley93}
S.~M. Buckley.
\newblock Estimates for operator norms on weighted spaces and reverse {J}ensen
  inequalities.
\newblock {\em Trans. Amer. Math. Soc.}, 340(1):253--272, 1993.

\bibitem{CV}
C.~Castaing and M.~Valadier.
\newblock {\em Convex Analysis and Measurable Multifunctions}, volume 580 of
  {\em Lecture Notes in Mathematics}.
\newblock Springer-Verlag, Berlin-New York, 1977.

\bibitem{MR684636}
M.~Christ and R.~Fefferman.
\newblock A note on weighted norm inequalities for the {H}ardy-{L}ittlewood
  maximal operator.
\newblock {\em Proc. Amer. Math. Soc.}, 87(3):447--448, 1983.

\bibitem{CG}
M.~Christ and M.~Goldberg.
\newblock Vector {$A_2$} weights and a {H}ardy-{L}ittlewood maximal function.
\newblock {\em Trans. Amer. Math. Soc.}, 353(5):1995--2002, 2001.

\bibitem{MR687639}
R.~Coifman, P.~W. Jones, and J.~L. Rubio~de Francia.
\newblock Constructive decomposition of {BMO} functions and factorization of
  {$A_{p}$} weights.
\newblock {\em Proc. Amer. Math. Soc.}, 87(4):675--676, 1983.

\bibitem{coifman-fefferman74}
R.~R. Coifman and C.~Fefferman.
\newblock Weighted norm inequalities for maximal functions and singular
  integrals.
\newblock {\em Studia Math.}, 51:241--250, 1974.

\bibitem{MR3521084}
J.. Conde-Alonso and G.~Rey.
\newblock A pointwise estimate for positive dyadic shifts and some
  applications.
\newblock {\em Math. Ann.}, 365(3-4):1111--1135, 2016.

\bibitem{dcu-paseky}
D.~Cruz-Uribe.
\newblock Extrapolation and factorization.
\newblock In J.~Lukes and L.~Pick, editors, {\em Function spaces, embeddings
  and extrapolation X, Paseky 2017}, pages 45--92. Matfyzpress, Charles
  University, 2017.

\bibitem{MR3803292}
D.~Cruz-Uribe, J.~Isralowitz, and K.~Moen.
\newblock Two weight bump conditions for matrix weights.
\newblock {\em Integral Equations Operator Theory}, 90(3):Art. 36, 31, 2018.

\bibitem{dcu-martell-perez}
D.~Cruz-Uribe, J.~M. Martell, and C.~P{\'e}rez.
\newblock Sharp weighted estimates for classical operators.
\newblock {\em Adv. Math.}, 229:408--441, 2011.

\bibitem{cruz-martell-perezBook}
D.~Cruz-Uribe, J.~M. Martell, and C.~P{\'e}rez.
\newblock {\em Weights, {E}xtrapolation and the {T}heory of {R}ubio de
  {F}rancia}, volume 215 of {\em Operator Theory: Advances and Applications}.
\newblock Birkh\"auser/Springer Basel AG, Basel, 2011.

\bibitem{cruz-uribe-perez00}
D.~Cruz-Uribe and C.~P{\'e}rez.
\newblock Two weight extrapolation via the maximal operator.
\newblock {\em J. Funct. Anal.}, 174(1):1--17, 2000.

\bibitem{MR3818613}
A.~Culiuc, F.~Di~Plinio, and Y.~Ou.
\newblock Uniform sparse domination of singular integrals via dyadic shifts.
\newblock {\em Math. Res. Lett.}, 25(1):21--42, 2018.

\bibitem{davey2022exponential}
B.~Davey and J.~Isralowitz.
\newblock Exponential decay estimates for fundamental matrices of generalized
  {S}chr\"odinger systems.
\newblock {\em preprint}, 2022.
\newblock arXiv:2207.05790.

\bibitem{MR4245601}
F.~Di~Plinio, T.~Hyt\"{o}nen, and K.~Li.
\newblock Sparse bounds for maximal rough singular integrals via the {F}ourier
  transform.
\newblock {\em Ann. Inst. Fourier (Grenoble)}, 70(5):1871--1902, 2020.

\bibitem{DPTV-2024}
K.~Domelevo, S.~Petermichl, S.~Treil, and A.~Volberg.
\newblock The matrix {$A_2$} conjecture fails, i.e., $3/2>1$.
\newblock {\em preprint}, 2024.
\newblock arXiv:2402.06961.

\bibitem{MR2140200}
O.~Dragi\v{c}evi\'{c}, L~Grafakos, M.~C. Pereyra, and S.~Petermichl.
\newblock Extrapolation and sharp norm estimates for classical operators on
  weighted {L}ebesgue spaces.
\newblock {\em Publ. Mat.}, 49(1):73--91, 2005.

\bibitem{duoandikoetxea01}
J.~Duoandikoetxea.
\newblock {\em Fourier analysis}, volume~29 of {\em Graduate Studies in
  Mathematics}.
\newblock American Mathematical Society, Providence, RI, 2001.

\bibitem{MR2104276}
M.~Frazier and S.~Roudenko.
\newblock Matrix-weighted {B}esov spaces and conditions of {$A_p$} type for
  {$0<p\leq1$}.
\newblock {\em Indiana Univ. Math. J.}, 53(5):1225--1254, 2004.

\bibitem{garcia-cuerva87}
J.~Garc{\'{\i}}a-Cuerva.
\newblock Jos\'e {L}uis {R}ubio de {F}rancia (1949--1988).
\newblock {\em Collect. Math.}, 38(1):3--15, 1987.

\bibitem{garcia-cuerva-rubiodefrancia85}
J.~Garc{\'{\i}}a-Cuerva and J.~L. Rubio~de Francia.
\newblock {\em Weighted norm inequalities and related topics}, volume 116 of
  {\em North-Holland Mathematics Studies}.
\newblock North-Holland Publishing Co., Amsterdam, 1985.

\bibitem{MR2015733}
M.~Goldberg.
\newblock Matrix {$A_p$} weights via maximal functions.
\newblock {\em Pacific J. Math.}, 211(2):201--220, 2003.

\bibitem{grafakos08b}
L.~Grafakos.
\newblock {\em Modern Fourier Analysis}, volume 250 of {\em Graduate Texts in
  Mathematics}.
\newblock Springer, New York, 2nd edition, 2008.

\bibitem{MR944321}
E.~Harboure, R.~A. Mac\'{\i}as, and C.~Segovia.
\newblock Extrapolation results for classes of weights.
\newblock {\em Amer. J. Math.}, 110(3):383--397, 1988.

\bibitem{MR4584703}
P.~Honz\'{\i}k.
\newblock An example of a singular integral and a weight.
\newblock {\em Int. Math. Res. Not. IMRN}, (9):7391--7398, 2023.

\bibitem{hytonenP2010}
T.~Hyt\"{o}nen.
\newblock The sharp weighted bound for general {C}alder\'{o}n-{Z}ygmund
  operators.
\newblock {\em Ann. of Math. (2)}, 175(3):1473--1506, 2012.

\bibitem{MR3204859}
T.~Hyt\"{o}nen.
\newblock The {$A_2$} theorem: remarks and complements.
\newblock In {\em Harmonic analysis and partial differential equations}, volume
  612 of {\em Contemp. Math.}, pages 91--106. Amer. Math. Soc., Providence, RI,
  2014.

\bibitem{hytonen-perez-treil-volbergP}
T.~Hy\"tonen, C.~P\'erez, S.~Treil, and A.~Volberg.
\newblock Sharp weighted estimates for dyadic shifts and the {$A_2$}
  conjecture.
\newblock {\em J. Reine Angew. Math.}, 687:43--86, 2014.

\bibitem{MR4030471}
J.~Isralowitz and K.~Moen.
\newblock Matrix weighted {P}oincar\'{e} inequalities and applications to
  degenerate elliptic systems.
\newblock {\em Indiana Univ. Math. J.}, 68(5):1327--1377, 2019.

\bibitem{Jones}
P.~W. Jones.
\newblock Factorization of {$A_{p}$} weights.
\newblock {\em Ann. of Math. (2)}, 111(3):511--530, 1980.

\bibitem{MR706075}
J.-L. Journ\'{e}.
\newblock {\em Calder\'{o}n-{Z}ygmund operators, pseudodifferential operators
  and the {C}auchy integral of {C}alder\'{o}n}, volume 994 of {\em Lecture
  Notes in Mathematics}.
\newblock Springer-Verlag, Berlin, 1983.

\bibitem{lacey-petermichl-reguera2010}
M.~Lacey, S.~Petermichl, and M.~Reguera.
\newblock Sharp ${A}_2$ inequality for {H}aar shift operators.
\newblock {\em Math. Ann.}, 348(1):127--141, 2010.

\bibitem{MR3127380}
A.~Lerner.
\newblock On an estimate of {C}alder\'{o}n-{Z}ygmund operators by dyadic
  positive operators.
\newblock {\em J. Anal. Math.}, 121:141--161, 2013.

\bibitem{lerner-IMRN2012}
A.~Lerner.
\newblock A simple proof of the {$A_2$} conjecture.
\newblock {\em Int. Math. Res. Not.}, 23(3):3159--3170, 2013.

\bibitem{MR4007575}
A.~Lerner and F.~Nazarov.
\newblock Intuitive dyadic calculus: the basics.
\newblock {\em Expo. Math.}, 37(3):225--265, 2019.

\bibitem{Lern2012}
A.~K. Lerner.
\newblock On an estimate of {C}alder\'{o}n-{Z}ygmund operators by dyadic
  positive operators.
\newblock {\em J. Anal. Math.}, 121:141--161, 2013.

\bibitem{muckenhoupt72}
B.~Muckenhoupt.
\newblock Weighted norm inequalities for the {H}ardy maximal function.
\newblock {\em Trans. Amer. Math. Soc.}, 165:207--226, 1972.

\bibitem{MR3689742}
F.~Nazarov, S.~Petermichl, S.~Treil, and A.~Volberg.
\newblock Convex body domination and weighted estimates with matrix weights.
\newblock {\em Adv. Math.}, 318:279--306, 2017.

\bibitem{NPST-2022}
F.~Nazarov, S.~Petermichl, K.~A. \v{S}kreb, and S.~Treil.
\newblock {The matrix-weighted dyadic convex body maximal operator is not
  bounded}.
\newblock {\em preprint}, 2022.
\newblock arXiv 2201.10337.

\bibitem{MR1428988}
F.~Nazarov and S.~Treil.
\newblock The hunt for a {B}ellman function: applications to estimates for
  singular integral operators and to other classical problems of harmonic
  analysis.
\newblock {\em Algebra i Analiz}, 8(5):32--162, 1996.

\bibitem{MR4000248}
Z.~Nieraeth.
\newblock Quantitative estimates and extrapolation for multilinear weight
  classes.
\newblock {\em Math. Ann.}, 375(1-2):453--507, 2019.

\bibitem{nieraeth-thesis}
Z.~Nieraeth.
\newblock {\em Sharp estimates and extrapolation for multilinear weight
  classes}.
\newblock PhD thesis, TU Delft, 2020.
\newblock https://doi.org/10.4233/uuid:192f633d-1bdb-440c-b435-7c0cd0d1f648.

\bibitem{MR2354322}
S.~Petermichl.
\newblock The sharp bound for the {H}ilbert transform on weighted {L}ebesgue
  spaces in terms of the classical {$A\sb p$} characteristic.
\newblock {\em Amer. J. Math.}, 129(5):1355--1375, 2007.

\bibitem{petermichl08}
S.~Petermichl.
\newblock The sharp weighted bound for the {R}iesz transforms.
\newblock {\em Proc. Amer. Math. Soc.}, 136(4):1237--1249, 2008.

\bibitem{petermichl-volberg02}
S.~Petermichl and A.~Volberg.
\newblock Heating of the {A}hlfors-{B}eurling operator: weakly quasiregular
  maps on the plane are quasiregular.
\newblock {\em Duke Math. J.}, 112(2):281--305, 2002.

\bibitem{MR3698161}
S.~Pott and A.~Stoica.
\newblock Bounds for {C}alder\'{o}n-{Z}ygmund operators with matrix {$A_2$}
  weights.
\newblock {\em Bull. Sci. Math.}, 141(6):584--614, 2017.

\bibitem{MR1928089}
S.~Roudenko.
\newblock Matrix-weighted {B}esov spaces.
\newblock {\em Trans. Amer. Math. Soc.}, 355(1):273--314 (electronic), 2003.

\bibitem{Rubio}
J.~L. Rubio~de Francia.
\newblock Factorization theory and {$A_{p}$} weights.
\newblock {\em Amer. J. Math.}, 106(3):533--547, 1984.

\bibitem{stein93}
E.~M. Stein.
\newblock {\em Harmonic analysis: real-variable methods, orthogonality, and
  oscillatory integrals}, volume~43 of {\em Princeton Mathematical Series}.
\newblock Princeton University Press, Princeton, NJ, 1993.
\newblock With the assistance of Timothy S. Murphy, Monographs in Harmonic
  Analysis, III.

\bibitem{MR1478786}
S.~Treil and A.~Volberg.
\newblock Continuous frame decomposition and a vector
  {H}unt-{M}uckenhoupt-{W}heeden theorem.
\newblock {\em Ark. Mat.}, 35(2):363--386, 1997.

\bibitem{MR1428818}
S.~Treil and A.~Volberg.
\newblock Wavelets and the angle between past and future.
\newblock {\em J. Funct. Anal.}, 143(2):269--308, 1997.

\bibitem{MR1423034}
A.~Volberg.
\newblock Matrix {$A_p$} weights via {$S$}-functions.
\newblock {\em J. Amer. Math. Soc.}, 10(2):445--466, 1997.

\end{thebibliography}

\end{document}